\documentstyle[psfig,12pt]{amsart}
\newcommand{\p}{\partial}
\newcommand{\sign}{\operatorname{sign}}
\newcommand{\rank}{\operatorname{rank}}
\newcommand{\ep}{\varepsilon}
\newcommand{\Th}{{\bold\Theta}^3_{\Bbb Z}}
\newcommand{\bz}{{\Bbb Z}}
\renewcommand{\fbox}{}

\title[On the homology cobordism group of homology 3-spheres] 
{On the homology cobordism group\\ of homology 3-spheres} 
\author{Nikolai Saveliev}
\address{Department of Mathematics \\ 
University of Michigan \\
Ann Arbor, MI 48109} 
\address{\it E-mail address:\ saveliev@@math.lsa.umich.edu}
\date{}

\begin{document}
\begin{abstract}

In this paper we present our results on the homology cobordism 
group $\Th$
of the oriented integral homology 3-spheres. We specially emphasize the role
played in the subject by the gauge theory including Floer homology and 
invariants by Donaldson and Seiberg -- Witten.
\end{abstract}
\vspace{5mm}

\maketitle
A closed oriented 3-manifold $\Sigma$ is said to be an integral homology 
sphere if it has the same integral homology as the 3-sphere $S^3$.
Two homology spheres $\Sigma_0$ and $\Sigma_1$ are homology cobordant,
if there is a smooth compact oriented 4-manifold $W$ with $\p W=
\Sigma_0\cup-\Sigma_1$ such that $H_*(W,\Sigma_0;\bz )= H_*(W,\Sigma_1;\bz )
=0$. The set of all homology cobordism classes forms an abelian group 
$\Th$ with the group operation defined by a connected sum. Here, the zero 
element is homology cobordism class of 3-sphere $S^3$, and the additive 
inverse is obtained by a reverse of orientation. 

The Rochlin invariant $\mu$ is an epimorphism $\mu :\Th\rightarrow \bz_2$, 
defined by the formula
$$\mu(\Sigma)={1\over 8}\ \sign (W)\ \mod \ 2,$$
where $W$ is any smooth simply connected parallelizable compact manifold 
with $\p W=\Sigma$. This invariant is well-defined due to well-known Rochlin 
theorem \cite {R} which states that the signature $\sign (V)$ of any smooth 
simply connected {\it closed} parallelizable manifold $V$ is divisible by 16.

We will focus our attention on the following problem concerning the group 
$\Th$ and the homomorphism $\mu$ : does there exist an element of order two 
in $\Th$ with non-trivial Rochlin invariant ?  This is one of R.~Kirby 
problems \cite{K}, Problem 4.4, and positive answer to this problem would
imply \cite{GS}, in particular, that all closed topological $n$-manifolds  are
simplicially triangulable if $n\ge 5$ ( this is not true in dimension 4 where 
a counterexample is due to A.~Casson and M.~Freedman, see e.g. \cite{AK} ).

Our approach to this problem is trying to lift the Rochlin homomorphism to
integers. If we succeded in doing this, then all elements of finite order 
would have to lie in the kernel of $\mu$, and the Kirby problem would 
have a negative solution. In fact, our goal is less ambitious -- we define 
a lift of $\mu$ on a certain subgroup of $\Th$. This implies that there 
are no solutions to the problem inside this subgroup. This is done in 
Section 1 with help of the $\bar\mu$--invariant introduced by W.~Neumann and
L.~Siebenmann in 1978, see \cite{N} and \cite{Sb}, for the so called plumbed 
homology spheres. 

Section 2 is devoted to an interesting observation about Floer homology of
homology 3-spheres homology cobordant to zero in the group $\Th$. We manage 
to prove a two-periodicity in Floer homology groups of such homology spheres 
in a series of examples. This leads us to defining a new integral invariant 
$\nu$ which vanishes for all these examples. Unfortunately, Floer homology is 
too complicated an object to deal with, so one still cannot say much more 
about the homology cobordism invariance of $\nu$ in general. 

In Section 3 we compare the invariants $\bar\mu$ and $\nu$ on the class of
plumbed homology 3-spheres and discuss possible links to branched coverings 
and Jones polynomial. 

\vspace{2mm}

Research at MSRI is supported in part by NSF grant DMS-9022140.

\vspace{5mm}

\section{Homology cobordisms of plumbed homology 3--spheres}
\vspace{4mm}

We first recall the definition of the invariant $\bar \mu$ by W. Neumann 
\cite{N}. Note that our definition differs from the original one by factor 
$1/8$.

Let $\Gamma$ be a plumbing graph, that is a (not necessarily connected) 
graph with no cycles, each of whose vertices carries an integer weight 
$e_i\ ,\ i=1,\ldots,s$. The matrix $A(\Gamma) = (a_{ij})_{i,j=1,\ldots,s}$ 
with the entries
$$
   a_{ij}  =  \cases e_i, &\text{ if }i=j \\
                     1, &\text{ if } i \text{ is connected to } j \text{ by an edge }\\
                     0, &\text{ otherwise }\endcases
$$
is the intersection matrix for 4-dimensional manifold $P(\Gamma)$, obtained
by plum\-bing $D^2$-bundles over 2-spheres according to $\Gamma$. 
This manifold is simply connected.

If $\Gamma$ is a plumbing graph as above, then $M(\Gamma)=\p P(\Gamma)$ is an 
integral homology sphere if and only if det$A(\Gamma)=\pm 1$. For example, 
all Seifert fibered homology spheres $\Sigma (a_1,\ldots,a_n)$ are of the 
form $\p P(\Gamma)$ where $\Gamma$ is a star-shaped graph, see \cite{NR}. 

If $M (\Gamma)$ is a homology sphere, there exists a unique homology class 
$w\in H_2 (P(\Gamma);\allowbreak \bz)$ satisfying the following two conditions. First, 
$w$ is {\it characteristic}, that is (dot represents intersection number)
\begin{equation}
 w . x \equiv x . x \quad \text{ mod\ 2 \quad for\ all }\ x\in
H_2 (P(\Gamma); \Bbb Z),
\end{equation}
\noindent and second, all coordinates of $w$ are either 0 or 1 in the natural
basis of $H_2 (P(\Gamma);\Bbb Z)$.  We call $w$ {\it the integral Wu class} for
$P (\Gamma)$.  Due to W.~Neumann \cite{N}, the integer $\text{sign} P(\Gamma) - 
w.w$ only depends on $M(\Gamma)$ and not on $\Gamma$. This integer is divisible
by 8, see \cite{Sr}, so one can define the Neumann-Siebenmann invariant by the
formula
\begin{equation}
\bar\mu(M(\Gamma)) = {1\over 8}(\text{ sign } P(\Gamma) - w.w).
\end{equation}
\noindent The class $w$ is spherical, so its modulo 2 reduction is the usual 
Rochlin invariant $\mu$, see \cite{NR}. 

W.~Neumann conjectured in \cite{N} that the number $\bar\mu$ is in fact a
homology cobordism invariant. This conjecture has been checked for all known 
examples of plumbed homology spheres homology cobordant to zero, including 
ones from the lists of A.~Casson and J.~Harer \cite{CH} and R.~Stern \cite{St}. 
The list in \cite{CH} includes Seifert fibered homology spheres 
$\Sigma(p,ps\pm 1,ps\pm 2),\ p\ge 3$ odd, $s\ge 1$, and $\Sigma(p,ps-1,ps+1),
\ p\ge 2$ even, $s\ge 1$ odd. All Stern's examples are of the form
$\Sigma(a_1, a_2, 2a_1a_2\pm a_3)$ for certain $\Sigma(a_1, a_2, a_3)$ in the 
Casson-Harer list. Now, gauge theory enables us to prove more general results
which weigh for positive answer to the conjecture.

\vspace{5mm}
\noindent{\bf Algebraic links.}
Our first result is concerned with the links of algebraic singularities. 
It was shown in \cite{EN} that a plumbed homology sphere $\Sigma$ is the link of 
an algebraic singularity if and only if there exists a plumbing graph 
$\Gamma$ such that the manifold $P(\Gamma)$ with $\Sigma = \p P(\Gamma)$ is 
negative definite. The simplest case is the Seifert fibered case: any Seifert 
fibered homology sphere $\Sigma (a_1,\ldots,a_n)$ is the link of the 
singularity of $f^{-1}(0)$ where $f:{\Bbb C}^n \to {\Bbb C}^{n-2}$ 
is a map of the form
$$f(z_1,\ldots,z_n) = \left( \sum_{k=1}^n b_{1,k} z_k^{a_k},\ldots,
                             \sum_{k=1}^n b_{n-2,k} z_k^{a_k}\right)
$$
with sufficiently general coefficient matrix $(b_{i,j})$, see \cite{NR}.
For instance,
$$\Sigma (p,q,r) = \{z\in {\Bbb C}^3 |\ \| z \|=\ep \text { and }
                                      z_1^p + z_2^q + z_3^r = 0\} $$
\noindent for $\ep>0$ small enough.

The next simplest case is the following: if $p,q,r$ are pairwise relatively
prime integers, as are $p',q',r'$, then the homology sphere $\Sigma$ obtained
by splicing $\Sigma (p,q,r)$ and $\Sigma (p',q',r')$ along the singular
fibers of degrees $r$ and $r'$, is the link of singularity if and only if
$rr' > pp'qq'$, see \cite{NW},\ \S\,4.

\vspace{4mm}
\noindent{\bf Theorem 1.} {\it Let a homology sphere $\Sigma$ be the link of 
an algebraic singularity. If $\Sigma$ is homology cobordant to zero then 
$\bar\mu (\Sigma) \ge 0$.}
\vspace{3mm}

\noindent{\it Proof.} Since $\Sigma$ is an algebraic link, one may assume that 
$\Sigma$ is the boundary of a plumbed negative definite 4-manifold 
$P (\Gamma)$. Suppose that $\Sigma$ bounds a smooth homology ball $M$. 
Let us consider the manifold $W = P(\Gamma) \cup_{\Sigma} (-M)$. This is a 
smooth closed oriented manifold whose intersection form is naturally 
isomorphic to the intersection form of $P (\Gamma)$, in particular, is 
negative definite. By S.~Donaldson's Theorem 1 from \cite{D}, this form is 
diagonalizable over integers.

We use this fact to evaluate $\bar\mu (\Sigma)$. Since $\sign P(\Gamma) =
-s$ where $s$ is the number of vertices in the graph $\Gamma$, we only need
to find the Wu class $w$. In the standard basis associated with plumbing, the 
matrix $A$ of the intersection form of $P (\Gamma)$ takes the form $A = U^t
(-E) U$ where $U\in \operatorname{SL}_s (\Bbb Z)$ and $E$ is the identity
matrix. The defining relation (1) translates to
$$w^t U^t (-E) U x \equiv x^t U^t (-E) U x \quad
            \text{ mod\ 2 \quad for\ all }\ x\in H_2 (P(\Gamma); \Bbb Z),$$
or, equivalently, 
$$(Uw)^t (-E) y \equiv y^t (-E) y \quad
            \text{ mod\ 2 \quad for\ all }\ y\in H_2 (P(\Gamma); \Bbb Z).$$
Therefore, $Uw$ is characteristic for $-E$, in particular, all coordinates of
$Uw$ are odd. Now, we have
$$w.w = w^t U^t(-E)U w = - (U w)^t(U w)$$
which is equal to the negative of the standard Euclidean square of
vector $U w$. Since all coordinates of $U w$ are odd, $w.w \le - s$, and
therefore $\bar\mu (\Sigma)\ge 0$.
\hfill\qed

\vspace{3mm}
\noindent{\bf Corollary.} {\it If plumbed $\Sigma$ is an algebraic link and 
$\bar\mu(\Sigma) < 0$, then $\Sigma$ has infinite order in the group
$\Th$.}
\hfill\qed

\vspace{3mm}
\noindent{\bf Example.} For any relatively prime integers $p, q > 0$, one 
can easily see that $\bar\mu (\, \Sigma(p, q,\allowbreak pq-1)) < 0$. Therefore, all 
these homology spheres are of infinite order in $\Th$. In particular, the
Poincar\'e homology sphere $\Sigma (2, 3, 5)$ is of infinite order.
\vspace{2mm}

\noindent{\bf Example.} The $\bar\mu$-invariant works not only for Seifert 
spheres, as it is seen from the following example. Let $\Sigma$ be a splice 
\cite{EN} of $\Sigma(4,7,9)$ and $\Sigma(2,3,25)$ along singular fibers of 
degrees 9 and 25. This manifold is an algebraic link. By using the additivity 
of $\bar\mu$ proven in \cite{S5}, we find that
$$\bar\mu(\Sigma) = \bar\mu(\Sigma(4,7,9)) + \bar\mu(\Sigma(2,3,25))
               = -2 + 0 < 0.$$
Therefore, $\Sigma$ has infinite order in $\Th$.
\vspace{5mm}

Many Seifert homology spheres $\Sigma (a_1,\ldots,a_n)$ having infinite order 
in the group $\Th$ by the corollary above, can also be detected by the 
Fintushel-Stern invariant
$$R(a_1,\ldots,a_n) = \frac 2{a} - 3 + n +
\sum_{i=1}^n \frac 2{a_i} \sum_{k=1}^{a_i -1} 
\cot \left(\frac {\pi a k}{a_i^2}\right)
\cot \left(\frac {\pi k}{a_i}\right) \sin^2 \left(\frac {\pi k}{a_i}\right) $$
where $a = a_1\cdot\ldots\cdot a_n$.
A theorem in \cite{FS2} says that if $\Sigma(a_1,\ldots,a_n)$ bounds a homology
ball then $R(a_1,\ldots,a_n) < 0$. There are however Seifert spheres which are
not homology cobordant to zero and which can be detected by $\bar\mu$ and not 
by $R$, and vice versa.

\vspace{3mm}
\noindent{\bf Example.} Both Seifert spheres $\Sigma (2,11,19)$ and 
$\Sigma (3,5,7)$ are not homology cobordant to zero. As to the former one, this
follows from the fact that $\bar\mu(\Sigma(2,11,19))=-1$ is negative ( though 
$R(2,11,19) = -1$); the latter one has $\bar\mu(\Sigma(3,5,7))=0$ and the
result follows from $ R(3,5,7) = 1$. 
\vspace{5mm}

Let $P(\Gamma)$ be a resolution of singularity whose link $\Sigma$ is a 
homology 3--sphere. Let $K$ be the Poincar\'e dual of the canonical class of 
a complex structure on $P(\Gamma)$, so it is characteristic. If $\Sigma$ is 
homology cobordant to zero then the same argument as in Theorem 1 shows that 
the number $\sign P - K . K$ should be non-negative (in fact, this is true
for {\it any} characteristic class $K$, not necessarily canonical). 

On the other hand, one can compute $\sign P - K . K$ easily in terms 
of the Milnor fiber $M$ of the singularity. Due to \cite{Mo}, 
$$(3\sign + 2e) (M) - K' . K' = (3\sign + 2e) (P) - K . K$$
where $e$ stands for Euler characteristic, and $K'$ is the canonical class of
$M$. Since manifold $P$ is negative definite, we get
$$\sign P - K . K = (3\sign + 2e) (M) - K' . K' - 2.$$
One can get more precise result for Seifert manifolds.

\vspace{4mm}
\noindent{\bf Theorem 2.} {\it If Seifert fibered homology sphere
$\Sigma (a_1,\ldots,a_n)$ is homology cobordant to zero then
$$\lambda (\Sigma (a_1,\ldots,a_n))\le \frac 1 {12} \prod_{i=1}^n (a_i - 1)$$
where $\lambda$ is the Casson invariant.}
\vspace{3mm}

\noindent{\it Proof.} In the special case of $\Sigma (a_1,\ldots,a_n)$, 
Milnor fiber $M$ is simply-connected and its canonical class $K'$ vanishes, 
see \cite{Dv}. In addition, 
$$\operatorname{rk} H_2 (M) = \prod_{i=1}^n (a_i - 1), \text{ and }
                       \sign M = -8 \lambda (\Sigma (a_1,\ldots,a_n)).$$
The first formula is due to J.~Milnor \cite{M}, and the second one to
Fintushel-Stern \cite{FS1} and Neumann-Wahl \cite{NW}. 
\hfill\qed

\vspace{3mm}
\noindent{\bf Example.} For any relatively prime integers $p,q$ and any 
positive integer $k$, $$\lambda (\Sigma (p,q,2pqk-1)) = k (p^2-1)(q^2-1)/12.$$
Therefore, homology sphere $\Sigma (p,q,2pqk-1)$ is not homology cobordant to 
zero due to Theorem 2 though its $\bar\mu$--invariant vanishes. In fact, this 
could be also proven with help of the Fintushel-Stern $R$-invariant. 

\vspace{5mm}
\noindent{\bf More constraints from even plumbing.}
Now, we have to deal with those plumbed homology spheres whose
$\bar\mu$--invariant is positive, and therefore, they are not prohibited from
bounding a homology ball by Theorem 1. A good example of this kind is 
$\Sigma (2,3,7)$ with $\bar\mu =1$. This homology sphere bounds the following
plumbing

\vspace{10mm}
\centerline{
\fbox{
\begin{picture}(200,50)
    \put(-5,5)        {\fbox{\psfig{figure=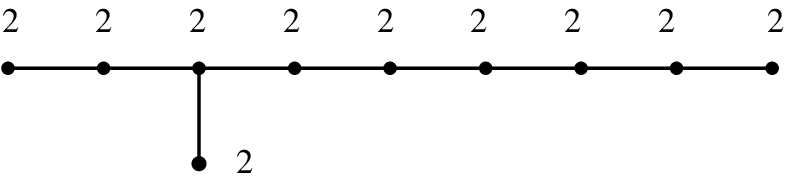}}}
\end{picture}
}
}
\vspace{3mm}

\noindent Its intersection form is even and indefinite, hence is isomorphic 
over integers to the form $E_8\oplus H$ where $E_8$ is the positive definite 
$E_8$--form, and $H$ is the hyperbolic $2\times 2$--matrix, see \cite{Sr}. 
Suppose that a double of $\Sigma (2,3,7)$ is homology 
cobordant to zero. Then, by gluing the homology ball and double of the plumbing
along $\Sigma (2,3,7)\# \Sigma (2,3,7)$, we get a closed smooth manifold with 
intersection form $2 E_8\oplus 2 H$ which is prohibited by S.~Donaldson's 
Theorem 2, see [D]. Therefore, $\Sigma (2,3,7)$ is not of the second order 
in $\Th$. 

The Seiberg-Witten gauge theory led M.~Furuta \cite{Fu}, see also \cite{A}, to 
a proof of the 10/8-conjecture saying that if a simply-connected closed smooth 
4-manifold has an even intersection form $Q$ then 
$$\frac {\rank (Q)}{|\sign (Q)|} > \frac {10}{8}.$$
One can use this result to show that in fact none of $\Sigma (2,3,7)$ multiples
can bound a contractible manifold. Because if any of them did, there would 
exist a closed smooth simply-connected manifold with rank/sign = 10/8.

One can use similar argument with even plumbings for many other Seifert fibered 
homology spheres with $\bar\mu$ positive.
Unfortunately, the construction does not work for {\it all} homology spheres
because the number of hyperbolics in even plumbing may be too large, and in 
order to get an even form with desirable quotient rank/sign one may need to 
surger some of them out (at the expense of losing the plumbing structure). In 
our paper [S2] we do this for Seifert fibered homology spheres 
$\Sigma (p,q,pqk\pm 1),\ k > 0$ odd, which can be thought of as manifolds 
obtained by $(-1/k)$--surgery along a $(p,q)$--torus knot. 

\vspace{5mm}
\noindent{\bf Theorem 3.} {\it Let $p, q > 1$ be relatively prime integers, 
and $k > 0$ an odd integer. Then neither $\Sigma (p,q,pqk+1)$ nor any of its
multiples can bound a smooth contractible 4--manifold.}
\vspace{5mm}

Note that $\bar\mu (\Sigma (p,q,pqk+1)) > 0$ for all $p,q,k$ as in the theorem,
and that a similar result for $\Sigma (p,q,pqk-1)$ followed from Theorem 1. 
Theorem 3 can not be generalized for $\Sigma (p,q,pqk+1)$ with $k$ even; for 
instance, $\Sigma (2,3,13)$ is known to bound a contractible manifold 
\cite{AK}. 

We refer for a complete proof of Theorem 3 to [S2], and only draw the final 
surgery diagrams in the simplest case of $p=2$. If $q\equiv 3\mod 4$, then 
$\Sigma (2,q,2qk+1)$ is the boundary of the following plumbing:

\vspace{10mm}
\centerline{
\fbox{
\begin{picture}(200,50)
    \put(5,0)        {\fbox{\psfig{figure=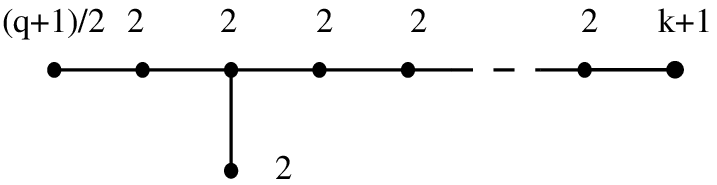}}}
\end{picture}
}
}
\vspace{3mm}

\noindent whose intersection form is isomorphic to $\frac {q+1}{4}\cdot E_8
\oplus H$. If $q\equiv 1\mod 4$, some hyperbolics have to be surgered out. 
The homology sphere $\Sigma (2,q,2qk+1)$ is then surgery on the link shown 
in Figure 1.

\vspace{10mm}
\centerline{
\fbox{
\begin{picture}(250,100)
    \put(25,0)       {\fbox{\psfig{figure=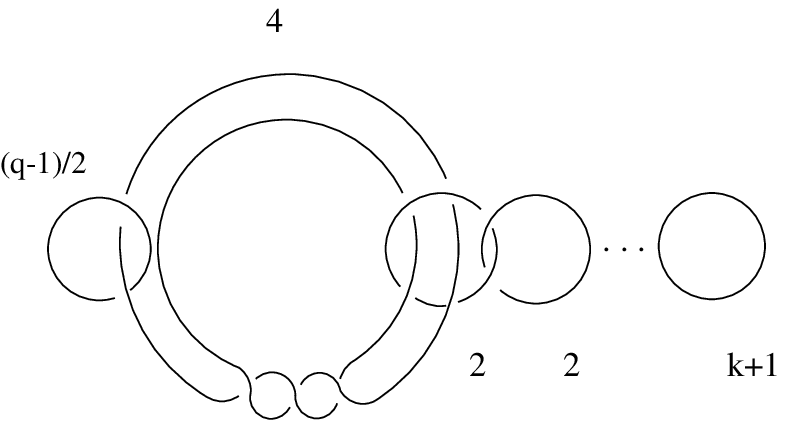}}}
    \put(120,-20)    {\makebox(0,0)[t]{\bf Figure 1}}
\end{picture}
}
}
\vspace{15mm}

\noindent The corresponding intersection form is isomorphic to $\frac
{q-1}{4}\cdot E_8\oplus H$.

\vspace{3mm}
\section{Floer homology and invariants of homology cobordism}
\vspace{2mm}

Floer homology groups $I_n(\Sigma),\ n = 0,\ldots,7$, are abelian groups, 
intrinsically associated to an oriented $\bz$-homology 3-sphere $\Sigma$. 
The definition of Floer homology in [F1] makes essential use of gauge 
theory on 3-- and 4--dimensional manifolds, so they are generally difficult 
to compute. 

We use the Floer exact triangle [F2] to compute $I_* (\Sigma)$ for several 
classes of integral homology spheres $\Sigma$ unified by one property -- 
they all are homology cobordant to zero. Among these homology spheres are: 
Mazur homology spheres \cite{AK}; Seifert fibered homology spheres which 
bound smooth contractible manifolds \cite{CH} and \cite{St}, homology 
spheres obtained by $(\pm 1)$-Dehn surgery on some slice knots. In short, 
the main results of our computations can be summarized as follows.  

\vspace{4mm}
\noindent{\bf Theorem 4.} {\it For all integral homology spheres $\Sigma$ from 
the list above, Floer homology groups are 2--periodic, that is, 
$I_j (\Sigma) = I_{j+2} (\Sigma)$ for any $j$.}
\vspace{3mm}

Below we will give more specific description of the homology spheres in 
Theorem 4, but first we would like to use this theorem to justify 
introduction of a new invariant of oriented integral homology 3--spheres.
For any such sphere $\Sigma$, we will define a number $\nu (\Sigma)$ 
keeping in mind the following conditions we would like it ideally to 
satisfy:

\quad {\it (1)} $\nu$ is a $\Bbb Z$-homology cobordism invariant,

\quad {\it (2)} $\nu$ is additive with respect to connected sums of
$\Bbb Z$-homology spheres, and

\quad {\it (3)} $\nu$ is a lifting of the Rochlin homomorphism $\mu$ , i.e.
$\nu(\Sigma)\equiv \mu(\Sigma) \mod 2$.

Recall that according to C.~Taubes \cite{T} Casson $\lambda$--invariant 
equals one-half of the Euler characteristic of Floer homology,
\begin{equation}
\lambda(\Sigma)={\frac {1}{2}}\sum_{n=0}^7 (-1)^n 
                           \operatorname{ rk}_{\Bbb Q} I_n (\Sigma).
\end{equation}
We define our invariant $\nu$ in a manner similar to (3), i.e. as an algebraic 
sum of the ranks of Floer homology groups. In doing so, we observe that there 
is essentially only one choice of signs, besides that in (3), that makes the 
invariant sensitive to orientations (as it is required by conditions {\it (1)} 
and {\it (2)} above). We use this choice to define our invariant $\nu$ as
follows:
\begin{equation}
\nu(\Sigma)={\frac{1}{2}}\sum_{n=0}^7 (-1)^{\frac {(n+1)(n+2)}{2}} 
                            \operatorname{ rk}_{\Bbb Q} I_n (\Sigma).
\end{equation}                            

\noindent This is an invariant of orientation preserving homeomorphism just
because so is Floer homology. The most important problem, therefore, is to
check its homology cobordism invariance, and Theorem 4 gives some evidence 
towards that.                           

\vspace{4mm}
\noindent{\bf Corollary.} {\it The $\nu$-invariant vanishes on all integral 
homology spheres homology cobordant to zero listed in Theorem 4.}
\vspace{3mm}

Concerning the conditions {\it (2)} and {\it (3)}, we mention that 
$\nu (\Sigma)$ is obviously equal to $\mu (\Sigma)$ modulo 2 for homology
spheres whose Floer homology is 4--periodic. In general, this is an open
problem. In \cite{S3} we use Fukaya's spectral sequence in Floer homology 
\cite{Fk} to show that $\nu (\Sigma)$ is additive with respect to connected 
sums under some restrictions on Floer homology of the factors.

\vspace{5mm}
\noindent{\bf Mazur homology spheres.} By a Mazur homology sphere we mean an
integral homology 3-spheres $\Sigma$ obtained by the following construction.
Let $W$ be a compact contractible 4-manifold obtained by attaching a 
two-handle to $S^1\times B^3$ along its boundary as in Figure 2.

\vspace{10mm}
\centerline{
\fbox{
\begin{picture}(200,110)
    \put(35,0)       {\fbox{\psfig{figure=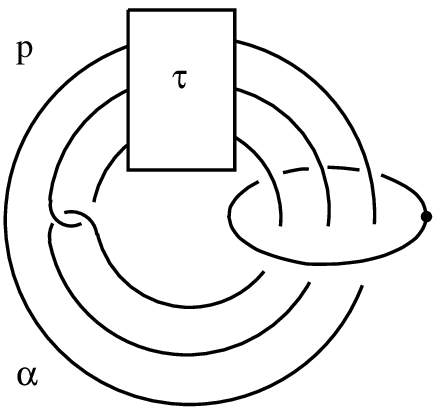}}}
    \put(95,-20)     {\makebox(0,0)[t]{\bf Figure 2}}
\end{picture}
\begin{picture}(200,110)
    \put(28,0)       {\fbox{\psfig{figure=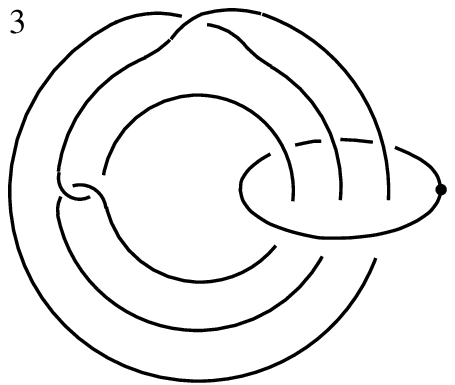}}}
    \put(95,-20)     {\makebox(0,0)[t]{\bf Figure 3}}
\end{picture}
}
}
\vspace{15mm} 

\noindent Here, the framing $p$ is an integer number, and $\tau$ is an
arbitrary tangle such that the component $\alpha$ of the link is connected. 
The manifold $W$ is a straight generalization of the well known Mazur manifold 
\cite{M} shown in Figure 3. The homology spheres $\Sigma$ in question are the 
boundaries of such 4-manifolds. For example, the boundary of the manifold
in Figure 3 is homeomorphic to Seifert fibered homology sphere 
$\Sigma (2,5,7)$, see \cite{AK}. The following result is proven in \cite{S3}.

\vspace{3mm}
\noindent {\bf Theorem 5.} {\it Let $\Sigma = \p W$ be the boundary of a
manifold $W$ in Figure 2. Then, for any tangle $\tau$ and for any $p$, Floer
homology groups $I_* (\Sigma)$ are two-periodic, that is $I_i (\Sigma) = 
I_{i+2} (\Sigma)$ for any $i$. Moreover, the groups $I_* (\Sigma)$ are
independent of the framing $p$.}
\vspace{3mm}

This theorem implies in particular that the $\nu$-invariant 
vanishes for all homology spheres as above. One can say more about 
$I_* (\Sigma)$ in the case when $\tau$ is a braid, see \cite{S2}.

\vspace{3mm}
\noindent {\bf Theorem 6.} {\it If $\tau$ is a braid in Figure 1, then for any
$p$, the group $I_i (\Sigma)$ is trivial if $i$ is odd, and $I_i (\Sigma) = 
\Bbb Z$ if $i$ is even, or vice versa, depending on the orientation.}

\vspace{5mm}
\noindent{\bf Casson-Harer and Stern homology spheres.} Mazur manifolds are 
examples of contractible manifolds having handlebody decomposition with one 
0-handle and one pair of 1- and 2-handles. In general, the geometric 
intersection of the 1- and 2-handles may be greater than or equal 3 (of course, 
the algebraic intersection still needs to be $\pm 1$). The Casson-Harer 
homology spheres \cite{CH} are of this sort. On the contrary, Stern homology 
spheres \cite{St} bound contractible manifolds with more than one pair of 1- 
and 2-handles. 

In \cite{S3}, we prove that Floer homology groups of all Seifert fibered 
homology spheres from \cite{CH} and \cite{St} known to bound contractible 
manifolds are 2--periodic. In our proof, we make use of Floer exact triangle 
applied to some special surgery on plumbed homology spheres that we introduce, 
and computations of Floer homology for Seifert manifolds in \cite{FS1}. Here 
are some explicit computations. 

\vspace{3mm}
\noindent {\bf Theorem 7.} {\it Let $p$ be any odd integer greater than 1, 
and $s\ge 1$ any integer. Then the Floer homology groups 
$I_j (\Sigma(p, ps\pm 1, ps\pm 2))$ are trivial for $j$ odd, and isomorphic 
to a free abelian group of rank $s(p^2-1)(ps\pm 3)/48$ for any even $j$. In 
particular, $\nu (\Sigma(p, ps\pm 1, ps\pm 2)) = 0$.}

\vspace{5mm}
\noindent{\bf Slice knots.} In \cite{S3} we also compute Floer homology
of homology spheres obtained by Dehn surgery on the knots $T(p,q)\# 
\allowbreak T(p,q)^*$
in $S^3$ where $T(p,q)$ is the left--handed $(p,q)$-torus knot, and 
$T(p,q)^*$ its mirror image. In the particular case of $(p,q)=(2,2n+1)$, the
Floer homology groups have been computed by A.~Stipsicz and Z.~Szab\'o in
\cite{SS} by different methods using the results of \cite{KKR} on splitting 
spectral flow.

It is worth mentioning that for any relatively prime integers $(p,q)$, the knot
$T(p,q)\# \allowbreak {T(p,q)}^*$ is slice, and that both $(+1)$-- and 
$(-1)$--surgeries
on it bound contractible manifolds. As it follows from Theorem 8 below, the
$\nu$-invariants vanish in both cases.

\vspace{3mm}
\noindent{\bf Theorem 8.} {\it Let $k=T(p,q)\# \allowbreak
T(p,q)^*$ be a connected sum of 
the $(p,q)$-torus knot $T(p,q)$ with its mirror image $T(p,q)^*$. Let
$M_n (p,q)$ be the homology sphere obtained by $(-1/n)$-surgery of $S^3$
along $k,\ n=1,2,\ldots$. Then $I_j (M_n (p,q))$ is trivial if $j$ is odd,
and isomorphic to a free abelian group of rank $n(p^2-1)(q^2-1)/24$ if $j$
is even.}
\vspace{3mm}

In our proof, we first use Floer exact triangle to discover some symmetries 
in Floer homology of amphicheiral knots, and then apply a description of the 
representation space of the fundamental group of a splice of Seifert fibered 
homology spheres along singular fibers, see \cite{S3}.

\vspace{3mm}
\section{Comparing the invariants}
\vspace{2mm}

Our computations of both $\nu$-- and $\bar\mu$--invariants suggest that these
invariants coincide for plumbed homology spheres. In comparing the two 
invariants, we use the well-known fact, see \cite{EN} and \cite{Sb}, that any 
plumbed homology sphere can be obtained by pasting together graph links which 
are Seifert fibered homology spheres with a collection of Seifert fibers. This 
operation is called splicing \cite{EN}. With respect to this operation, the 
$\bar \mu$-invariant is additive, see e.g. [S5]. In particular, if 
$\Sigma(a_1,\ldots,a_n)$ is a Seifert fibered homology sphere with Seifert 
invariants $a_1,\ldots,a_n$ and, for any $2\le j \le n-2$, the integers 
$q=a_1\cdots a_j$ and $p=a_{j+1}\cdots a_n$ are the products of the first $j$ 
and the last $(n-j)$ Seifert invariants, respectively, then

$$\bar \mu(\Sigma(a_1,\ldots,a_n))= \bar \mu(\Sigma(a_1,\ldots,a_j,p)) + 
                               \bar\mu(\Sigma(q,a_{j+1},\ldots,a_n)),$$
$$\nu(\Sigma(a_1,\ldots,a_n))= \nu(\Sigma(a_1,\ldots,a_j,p)) + 
                               \nu(\Sigma(q,a_{j+1},\ldots,a_n)).$$

\noindent The formula for $\nu$ follows from [S1]. This implies that the problem of 
identifying $\bar\mu$ and $\nu$ reduces, at least on the class of Seifert 
fibered homology spheres, to the identification problem for Seifert spheres 
with just three singular fibers.

One can also compare the behaviour of the two invariants with respect to Dehn 
surgery on a singular fiber in a Seifert fibered homology sphere. Let 
$\Sigma(a_1,\ldots,a_n)$ be a Seifert fibered homology sphere with Seifert 
invariants $a_1,\ldots,a_n$, and $a=a_1\cdots a_{n-1}$ be the product of the 
first $(n-1)$ Seifert invariants. Then

$$\bar \mu(\Sigma(a_1,\ldots,a_n))=\ \pm\ \bar \mu(\Sigma(a_1,\ldots,a_{n-1},
2a\pm a_n)).$$
$$\nu(\Sigma(a_1,\ldots,a_n))=\ \pm\ \nu(\Sigma(a_1,\ldots,a_{n-1},2a\pm
a_n)),$$

\noindent The first part of this result is proven in \cite{N}, the second in 
\cite{S3}. As a corollary, both invariants $\bar \mu$ and $\nu$ vanish on 
infinite series of Seifert fibered homology spheres $\Sigma(a_1,\ldots,a_{n-1},
 2ak\pm 1),\ k\in \bz$, where $a=a_1\cdots a_{n-1}$. 

Thus, in order to identify the two invariants for all Seifert fibered homology
spheres, one only needs to do this for homology spheres of the form
$\Sigma (p,q,r)$ for all pairwise relatively prime integers $p, q, r$ such that
$\max (p,q) < r < pq$. For example, all Casson-Harer homology spheres are of
this sort. The invariants $\bar\mu$ and $\nu$ for these manifolds coincide
(since they vanish). Our numerous attempts at proving this fact in general have
led us to a few more special cases when it holds, as well as to some its 
reformulations. One of them seems to be worth stating.

It is well-known that a Seifert fibered homology sphere 
$\Sigma (a_1,\ldots,a_n)$ is a doble cover over $S^3$ branched over Montesinos 
knot $K (a_1,\ldots,a_n)$, see e.g. \cite{BZ}. The $\bar\mu$--invariant can be
defined in terms of the signature of this knot, as it was originally done by
L.~Siebenmann \cite{Sb}. On the other hand, there exists a formula for Casson 
invariant of a double branched cover by D.~Mullins \cite{Mu}. These two facts
together with formulae (3) and (4) enable us to reformulate our conjecture
about the coincidence of $\bar\mu$ and $\nu$ for Seifert fibered homology 
spheres in the following equivalent form

$$\rank_{\Bbb Q} I_0 (\Sigma (a_1,\ldots,a_n)) = \pm 
\frac 1 {12} \cdot \frac d {dt} \ \ln V_K (-1)$$

\noindent where $V_K$ is the Jones polynomial of the Montesinos knot $K =
K(a_1,\ldots,a_n)$.

\end{document}